\documentclass{amsart}

\input{diagrams}

\usepackage{amssymb}
\usepackage{amsmath}
\usepackage{amsthm}
\usepackage{amscd}

\newcommand{\Z}{{\mathbb Z}}            
\newcommand{\Q}{{\mathbb Q}}            
\newcommand{\R}{{\mathbb R}}            
\newcommand{\C}{{\mathbb C}}            
\newcommand{\Zp}{{\mathbb Z_p}}         
\newcommand{\Qp}{{\mathbb Q_p}}         
\newcommand{\Fpbar}{{\overline{\F}_p}}           
\newcommand{\Kbar}{{\bar K}}            
\newcommand{\Kvbar}{{\bar K_v}}            
\newcommand{\Qbar}{{\overline{\Q}}}
\DeclareMathOperator{\unr}{{\rm nr}}    
\newcommand{\Qpunr}{\Q^{\unr}_p}        
\newcommand{\p}{{\mathfrak{p}}}         

\def\isomap{{\buildrel \sim\over\longrightarrow}} 

\newcommand{\E}{{\mathcal E}}           
\newcommand{\A}{{\mathcal A}}           
\newcommand{\Et}{\tilde{E}}             

\newcommand{\M}{{\mathcal M}}           

\renewcommand{\l}{\lambda}              
\newcommand{\s}{\sigma}                 

\newcommand{\jinv}{j_E^{-1}}            
\newcommand{\Gm}{{\mathbf G}_m}         
\newcommand{\G}{{\mathbf G}}    
\newcommand{\Gdual}{G^\vee}    
\newcommand{\Adual}{A^\vee}    
\newcommand{\Gmhat}{\widehat{{\mathbf G}}_m}    
\newcommand{\Ehat}{\widehat{\E}}                
\newcommand{\Ahat}{\widehat{\A}}                
\newcommand{\B}{\mathbf B}              

\DeclareMathOperator{\End}{{\rm End}}                   
\DeclareMathOperator{\Gal}{{\rm Gal}}           
\DeclareMathOperator{\tor}{{\rm tor}}           
\DeclareMathOperator{\ord}{{\rm ord}}           

\newcommand{\F}{{\mathbf F}}                    

\newcommand{\f}{\mathfrak{f}}           

\newcommand{\hhat}{\hat{h}}             
\newcommand{\Kab}{K^{\rm ab}}           

\newtheorem{theorem}{Theorem}[section]
\newtheorem{conj}[theorem]{Conjecture}
\newtheorem{cor}[theorem]{Corollary}
\newtheorem{lemma}[theorem]{Lemma}

\def\pf{{\sc Proof. }}
\def\qed{\hfill $\blacksquare$\smallskip}

\theoremstyle{remark}
\newtheorem{remark}[theorem]{Remark}

\newcounter{listcounter1}
\newenvironment{list1}{
  \begin{list}{\arabic{listcounter1}.\hfill}{
    \usecounter{listcounter1}
    \setlength{\leftmargin}{18pt}
    \setlength{\labelwidth}{18pt}
    \setlength{\labelsep}{0pt}
    \setlength{\topsep}{0pt}
  }
}{
  \end{list}
}

\begin{document}

\title[Canonical heights on elliptic curves over abelian extensions]
{Lower bounds for the canonical height on elliptic curves over abelian
extensions}


\date{Revised on December 16, 2002}

\author[Matthew Baker]{Matthew H. Baker}
\email{mbaker@math.uga.edu}

\address{Department of Mathematics, 
         University of Georgia, Athens, GA 30602-7403}

\thanks{The author's research was supported in part 
by an NSF Postdoctoral Research Fellowship.}

\begin{abstract}
Let $K$ be a number field and let $E/K$ be an elliptic curve.
If $E$ has complex multiplication, we show that
there is a positive lower bound for the canonical height of 
non-torsion points on $E$ defined over the maximal abelian extension 
$\Kab$ of $K$.
This is analogous to results of Amoroso-Dvornicich and 
Amoroso-Zannier for the multiplicative group.  We also show that if
$E$ has non-integral $j$-invariant (so that in particular $E$ does not have
complex multiplication), then there exists $C > 0$ such that there are only finitely many
points $P \in E(\Kab)$ of canonical height less than $C$.  This 
strengthens a result of Hindry and Silverman.
\end{abstract}

\maketitle

\section{Introduction}
\label{introsection}

Let $K$ be a number field, and let $\Kab$ denote the 
maximal abelian extension of $K$.  
Let $E/K$ be an elliptic curve, and let $\hhat : E(\Kbar) \to \R$
be the usual canonical height function
on $E(\Kbar)$.

Suppose that $E$ has complex multiplication by an order 
in the imaginary quadratic field $k$.  
If $K$ is any number field containing $k$,
the theory of complex multiplication tells us that all of the torsion 
points of $E(\Kbar)$ 
(which are exactly the points with canonical height zero) 
are defined over the maximal abelian extension $\Kab$ of $K$. 
By analogy with recent results of Amoroso--Dvornicich and Amoroso--Zannier
for the multiplicative group (see Theorem~\ref{ADtheorem} below), 
one might ask whether or not the canonical height function, when restricted
to $E(\Kab)$, takes arbitrarily small positive real values.  The 
following theorem answers this question.

\begin{theorem}
\label{maintheorem}
Let $K$ be a number field, and let $E/K$ be an elliptic curve having complex
multiplication.
Then there exists an effectively computable real constant $C>0$ (depending on $E/K$) 
such that $\hhat(P) \geq C$ for all non-torsion points $P\in E(\Kab)$.
\end{theorem}


\begin{remark}
We will prove a somewhat more precise statement (Theorem~\ref{precise}) in 
Section~\ref{proofsection}.
\end{remark}

\medskip

Theorem~\ref{maintheorem} is analogous to the following result for heights 
on the multiplicative group $\Gm$, which (although not stated in
precisely this form) is proved in \cite{AD} and \cite{AZ}:

\begin{theorem}
\label{ADtheorem}
Let $h$ denote the usual absolute logarithmic Weil height on $\Qbar$.  
Let $K$ be a number field, and let $\Kab$ be the maximal abelian extension of
$K$.  Then there exists an effectively computable constant $C>0$ 
(depending only on $K$) such that 
if $\alpha \in \Kab$ is not equal to zero or a root of unity, then
$h(\alpha) \geq C$.
\end{theorem}

\medskip

Now suppose that $E$ does not have complex multiplication.  
As a complement to Theorem~\ref{maintheorem}, we prove the following
result:

\begin{theorem}
\label{maintheorem2}
Suppose $K$ is a number field and $E$ is an elliptic curve with non-integral
$j$-invariant.  Then there exists a constant $C>0$ such that 
$\hhat(P) \geq C$ for all but finitely many points $P\in E(\Kab)$.
\end{theorem}

\begin{remark}
The hypothesis that $j_E$ is non-integral implies that $E$ does not have
complex multiplication.  The elliptic curves not covered by either
Theorem~\ref{maintheorem} or Theorem~\ref{maintheorem2} are exactly
those which have everywhere potential good reduction but do not
have complex multiplication.
\end{remark}

\begin{remark}
Theorem~\ref{maintheorem2} generalizes a result of Hindry-Silverman
(Corollary~0.3 of \cite{HS2}), which says that under the hypotheses of
Theorem~\ref{maintheorem2}, there exist effectively computable
constants $C,D>0$ such that any point $P \in E(\Kab)$ with $\deg(P) \geq
D$ satisfies the inequality
\[
\hhat(P) \geq \frac{C}{\deg(P)^{2/3}}.
\]

Silverman also proved in \cite{SilvermanDMJ} that if $E/K$ is an
elliptic curve without complex multiplication (over $\Kbar$), then
there exists an effectively computable constant $C>0$ such that every non-torsion point $P \in
E(\Kab)$ satisfies 
\[
\hhat(P) \geq \frac{C}{\deg(P)^{2}}.
\]

\end{remark}

\begin{remark}
The proof of Theorem~\ref{maintheorem2} does not provide an effective
algorithm for computing the constant $C$ in the statement of the theorem.
\end{remark}

\medskip

We conclude this section by placing Theorems~\ref{maintheorem} and 
\ref{maintheorem2} in a more general context.

Let $G$ be a variety defined over a field $K$, and let $L/K$ be
a field extension.
Suppose that $\hhat : G(L) \to \R$ is a function taking nonnegative
values and that $S \subseteq G(L)$.
We say that the pair $(S,\hhat)$ (or the set $S$, when the function
$\hhat$ is understood) has the {\em discreteness property} if there exists a
constant $C>0$ such that $\hhat(P) \geq C$ whenever $P \in S$ and
$\hhat(P) > 0$.
We say that the pair $(S,\hhat)$ has the {\em strong discreteness property} if 
in fact we have
\[
\liminf_{P \in S} \hhat(P) > 0,
\]  
i.e., if there exists a constant $C>0$ such that $\hhat(P) \geq C$ for
all but finitely many $P \in S$.


An example to keep in mind is the case where $G=A$ is an abelian variety defined over a
number field $K$ and
$\hhat : A(\Kbar) \to \R$ is the canonical height function corresponding
to some symmetric ample line bundle on $A$.
It is well known (see \cite[Section 4]{Poonen}, for example) that
whether or not $(S,\hhat)$ has the discreteness property (resp. the strong
discreteness property) is independent of which canonical height function one chooses.  
Another example is where $G=\Gm^n$ and 
$\hhat : G(\Qbar) \to \R$ is the function given 
$\hhat(x_1,\ldots,x_n) = h(x_1) + \cdots + h(x_n)$, where $h : \Qbar
\to \R$ is the usual logarithmic Weil height.

\medskip

Theorem~\ref{ADtheorem} implies that $\Gm^n(\Kab)$ has the 
discreteness property for all number fields $K$.
Theorem~\ref{maintheorem} can be rephrased as saying 
that if $E/K$ is an elliptic curve with complex multiplication, 
then $E(\Kab)$ has the discreteness property.  
Similarly, Theorem~\ref{maintheorem2} says that if $E/K$ is an elliptic
curve with non-integral $j$-invariant,
then $E(\Kab)$ has the strong discreteness property.

\medskip

We mention some known results related to Theorems~\ref{maintheorem} and
\ref{maintheorem2}.  Concerning torsion points of
abelian varieties defined over abelian extensions, J.~Zarhin \cite{Zarhin}
(see also Ruppert \cite{Ruppert}) has proved the following, using deep finiteness 
results of Faltings (\cite{Faltings}):

\begin{theorem}
\label{Ruppert-Zarhin}
Let $A$ be an abelian variety defined over a number field $K$.  
Then the torsion subgroup $A_{\tor} (\Kab )$ of $A(\Kab)$ 
is finite iff $A$ has no abelian subvariety
with complex multiplication over $K$.
\end{theorem}


There is also the following result due to S.~Zhang, which is a
consequence (using Weil restriction of scalars)
of the equidistribution theorems proved in \cite{Zhang} and \cite{SUZ}:

\begin{theorem}[Zhang]
Let $K$ be a number field, and let $A/K$ be an abelian variety.  
If $L \supseteq K$ and $L$ is a finite extension of 
a totally real field, 
then $A(L)$ has the strong discreteness property.
\end{theorem}


Since the maximal cyclotomic extension $K^{\rm cycl} := K(\mu_{\infty})$ 
is a finite extension of the totally real field $\Q(\mu_{\infty})^+$, 
Zhang's result implies in particular:

\begin{cor}
\label{Zhangcor}
Let $A/K$ be an abelian variety.  Then $A(K^{\rm cycl})$ has the strong
discreteness property.
\end{cor}

The weaker fact that the torsion subgroup 
$A_{\tor}(K^{\rm cycl})$ is finite is originally due to Ribet 
\cite{Ribet}.

In light of the above-mentioned results, we propose the following 
conjecture.  

\begin{conj}
\label{conjecture1}
Let $K$ be a number field, and let $A/K$ be an abelian variety.  Then:
\begin{list1}
\item $A(\Kab)$ has the discreteness property.
\item If $A$ does not contain an abelian subvariety having complex multiplication over $K$, 
then $A(\Kab)$ has the strong discreteness property.
\end{list1}
\end{conj}

Note that Theorems~\ref{maintheorem} and \ref{maintheorem2} 
are special cases of part (1) of this conjecture.
Together with Corollary~\ref{Zhangcor},
these results imply that Conjecture~\ref{conjecture1} is true when $A=E$ is
an elliptic curve and any one of the following three hypotheses is 
satisfied:
\begin{itemize}
\item $E$ has complex multiplication
\item $j_E$ is non-integral
\item $K=\Q$.
\end{itemize}

\medskip

We note that if $P_0 \in A(K)$ is a nontorsion point and $P_n \in
A(\Kbar)$ satisfies $2^n P_n = P_0$, then
$\hhat(P_n) = \frac{1}{4^n} \hhat(P_0) \to 0$ as
$n\to\infty$.  This is the simplest way to construct a sequence of
nontorsion points in $A(\Kbar)$ whose heights tend toward zero.
Conjecture~\ref{conjecture1} is therefore related to the principle that
when $n$ is large, the points $P_n$ tend to be defined over highly 
non-abelian extensions of $K$.

\medskip

To conclude this introduction, we pose the following question (which one could also ask 
in the context of semiabelian varieties):

\medskip

{\bf Question:} If $A/K$ is an abelian variety and 
$L$ is the extension of $K$ generated by all torsion points of $A(\Kbar)$,
is it true that $A(L)$ has the discreteness property?

\medskip


\section{Background facts on canonical heights}
\label{heightreviewsection}

If $K$ is a number field, we let $M_K$ denote the set of places (i.e.,
equivalence classes of absolute values) of $K$.  We will denote the archimedean
places by $M_K^\infty$ and the non-archimedean ones by $M_K^0$.  
If $v \in M_K$, we let $|\cdot |_v$ be the unique absolute value in the equivalence class $v$ which
extends the usual absolute value on $\Q_v$.
We also let $n_v = [K_v : \Q_v]$ denote the local degree at
$v$, and we let $d_v = \frac{[K_v : \Q_v]}{[K : \Q]}$ denote the local
degree divided by the global degree.

We have:
\begin{lemma}
\label{absolutevaluelemma}
If $L$ is a finite extension of the number field $K$ and $v$ is a fixed
place of $K$, then
\[
\sum_{w \in M_L, w \mid v} d_w = d_v.
\]
\end{lemma}

\pf Using the well-known fact that $\sum_{w \mid v} [L_w : K_v] = [L:K]$,
we have
\[
\sum_{w \mid v} \frac{[L_w : \Q_w]}{[L:\Q]}
= \sum_{w \mid v} \frac{[L_w : K_v][K_v : \Q_v]}{[L:\Q]}
= \frac{[K_v : \Q_v]}{[K:\Q]} 
\]
as desired.
\qed

\medskip

Let $E / K$ be an elliptic curve.  It is well known (see
  \cite[Chapter VI]{Silverman2} for definitions and further
  properties) that the canonical height function $\hhat : E(\Kbar) \to \R$ has
  the following properties:
\begin{itemize}
\item $\hhat(nP) = n^2 \hhat(P)$ for all $P \in E(\Kbar)$ and all
  integers $n$.
\item $\hhat(P) \geq 0$ for all $P \in E(\Kbar)$, and $\hhat(P) = 0$ if and only if $P$ is a torsion point.
\item (Northcott's theorem) Given any constant $M > 0$ and a natural
  number $d$,
  the set of points $P \in E(\Kbar)$ such that $\hhat(P) \leq M$ and
  $[K(P):K] \leq d$ is finite.
\item (Parallelogram law) For all $P,Q \in E(\Kbar)$,
\[
\hhat(P + Q) + \hhat(P - Q) = 2\hhat(P) + 2\hhat(Q).
\]
\end{itemize}

The canonical height $\hhat$ has a decomposition into N{\'e}ron local
height functions $\lambda_v$, one for each $v \in M_K$.  More precisely, for all $P
\in E(K), P \neq 0$ we have
\[
\hhat(P) = \sum_{v \in M_K} n_v \lambda_v (P),
\]
where $\lambda_v : E(K_v) - \{ 0 \} \to \R$ is defined and normalized as in
\cite[Chapter VI]{Silverman2}.  Note that in general, the N{\'e}ron local height
functions can take negative values, although if $v$ is a place of good
reduction for $E$, then $\lambda_v(P)$ is always nonnegative.  In
fact, suppose $E$ has good reduction at $v$ and that $\E$ is a minimal
Weierstrass model over the ring of integers $R_v$ in the completion
$K_v$.  Then if $v$ corresponds to the prime ideal $\p = (\pi)$ of $R_v$
and $P \in E(K_v) - \{ 0 \}$, we have 
\[
\lambda_v(P) = \log |\pi|_v^{m(P)},
\]
where $m(P)$ is the largest integer $m\geq 0$ such that $P \equiv 0$ (mod $\pi^m$).


The {N\'e}ron local heights behave functorially with respect to finite
extensions of $K$: if
$L/K$ is a finite extension and $P \in E(K)$, then
\[
\l_v(P) = \sum_{w | v} \frac{[L_w : K_v]}{[L:K]} \l_w(P).
\]

In particular, for each $v \in M_K$, $\lambda_v$ extends to a N{\'e}ron local height function
$\l_v : E(\overline{K_v}) - \{ 0 \} \to \R$.

\medskip

In order to prove Theorems~\ref{maintheorem} and 
Theorem~\ref{maintheorem2}, we will need known lower bounds
for certain average values of N{\'e}ron local height functions.
We consider the archimedean case first.  

\begin{lemma}
\label{Elkies}
Let $E/\C$ be an elliptic curve with $j$-invariant $j_E$, and let
$\l : E(\C) - \{ 0 \} \to \R$
denote the associated N{\'e}ron local height function.  
Let $P_1,\ldots,P_N$ be distinct points of $E(\C)$.  Then 
\[
\liminf_{N\to\infty} \frac{1}{N(N-1)} \sum_{i\neq j} \l(P_i - P_j)
\geq 0.
\]
More explicitly, we have
\[
\sum_{i\neq j} \l(P_i - P_j) \geq 
-\frac{1}{12}N\log N - \frac{1}{12}N\log^+|j_E| - \frac{11}{3}.
\]
\end{lemma}

\pf This follows from a result of Elkies which is
proved in a slightly different form in \cite[Theorem~VI.5.1]{Lang}.  The
necessary modifications to get the form quoted here are discussed in
\cite[Lemma 3]{HS3}.
\qed


It is clear from the definition of $\l$ that there exists a real constant $C_1$ such that $\l(P)
\geq -C_1$ for all $P \in E(\C) - \{ 0 \}$.
Taking $N=2$ in the above lemma, we obtain the
following explicit admissible value for $C_1$:

\begin{cor}
\label{Elkiescor}
Let $E/\C$, $\l$ be as above.
Then if $P\in E(\C) - \{ 0 \}$, we have
\[
\l(P) \geq -C_1,
\]
where $C_1 = \frac{1}{6}(C_2 + \log^+{|j_E|})$ and 
$C_2$ is a universal constant which can be taken to be
$\log 2 + \frac{22}{3}$.
\end{cor}


We now turn to the special case where $v$ is a non-archimedean place
of $K$ at which $E$ has split multiplicative reduction, so that in particular
the $j$-invariant of $E$ is non-integral at $v$.
Let $\nu = \ord_{v}(\jinv) > 0$, and note also for future reference that $\log |j_E|_v > 0$.

Let $\E$ be the N{\'e}ron model of $E$ over the ring of integers of $K_v$,
and let $\Et$ be its special fiber.
We denote by $\Phi$ the
group of connected components of $\Et$ and we fix an isomorphism
\[
\Phi \; \isomap \; \frac{(1/\nu)\Z}{\Z}
\]
(see \cite[Corollary IV.9.2]{Silverman2}).

Let $\l_v : E(K_v) - \{ 0 \} \to \R$ be the 
N{\'e}ron local height function at $v$.
Tate's explicit formula for $\l_v$ implies that if 
$r_v$ denotes the natural reduction 
homomorphism
\[
r_v : E(K_v) \to \Phi \; \isomap \; \frac{(1/\nu)\Z}{\Z},
\]
then 
\[
\l_v(P) \geq \frac{1}{2} \B_2(r_v(P))\log |j_E|_v
\]
for all nonzero points $P \in E(K_v)$.
Here $\B_2(t)$ is the periodic second Bernoulli polynomial
$\{ t \}^2 - \{ t \} + \frac{1}{6}$.  Note in particular that if $P
\in E(K_v) - \{ 0 \}$ reduces to the identity component of $\Et$ then 
\[
\l_v(P) \geq \frac{1}{12} \log |j_E|_v > 0.
\]

Using this description and a Fourier averaging technique, one can
prove the following analogue of Lemma~\ref{Elkies} (see \cite[Proposition~1.2]{HS2}):

\begin{lemma}
\label{HSlemma}
Let $v$ be a finite place of $K$ such that $j_E$ is non-integral at
$v$, and let $P_1,\ldots,P_N \in E(K_v)$ be distinct points.  
Let $\nu = \ord_v(\jinv)$.  Then
\[
\liminf_{N\to\infty} \frac{1}{N(N-1)} \sum_{i\neq j} \l_v(P_i - P_j)
\geq 0.
\]
More explicitly, we have
\[
\sum_{i\neq j} \lambda_v (P_i - P_j) \geq 
\frac{1}{12} \log |j_E|_v \left( (\frac{N}{\nu})^2 - N \right).
\]
\end{lemma}

\section{Preliminary Lemmas}
\label{lemmasection}

In this section, we prove some lemmas which will be needed for the
proof of Theorem~\ref{maintheorem}.
The first lemma is well-known, and concerns the formal group of an
abelian variety with ordinary reduction. 
For lack of a suitable reference, we sketch a proof.

\begin{lemma}
\label{formalgrouplemma}
Let $R$ be a strictly henselian discrete valuation ring 
with fraction field $K$ of characteristic zero and residue field $k$ of 
characteristic $p>0$.  
Let $A / K$ be an abelian variety with good ordinary reduction, 
and let $\A / R$ be its N{\'e}ron model.  
Let $\Ahat$ be the formal group 
(i.e., the formal completion along the identity section) of $\A$.
Then $\Ahat$ is isomorphic over $R$ to $(\Gmhat)^g$, where $g = \dim
A$ and $\Gmhat$ is the formal group of $\G_m$.
\end{lemma}

\pf (Sketch) Let $G / R$ be the 
$p$-divisible group of $\A$, so that by the Serre-Tate correspondence,
the formal group $\Ahat$ is the connected component $G^0$ of $G$.
Since $A$ has good ordinary reduction and $k$ is separably closed, $G$
contains a constant {\'e}tale $p$-divisible group isomorphic to $(\Q_p
/ \Z_p)^g$.

On the other hand, let $\Gdual$ be the $p$-divisible group of the dual abelian variety $\Adual$.
Then $G$ and $\Gdual$ are canonically Cartier dual to one another (see
\cite{Oda}).
Additionally, $A$ and $\Adual$ are isogenous over $\Kbar$, so
that $\Adual$ also has good ordinary reduction.  Therefore $\Gdual$
contains a copy of $(\Q_p / \Z_p)^g$,
and by Cartier duality, it follows that $G$ contains a split multiplicative
$p$-divisible group isomorphic to $(\G_m)^g$.

It follows by dimension considerations that $G \cong (\Q_p / \Z_p)^g
\times (\G_m)^g$ and therefore $\Ahat \cong (\Gmhat)^g$ as desired.

\qed

\medskip

Our proof of Theorem~\ref{maintheorem} will follow the same basic outline as
the proof of \cite[Theorem 1]{AD}, replacing the map
$x \mapsto x^p$ on $\Gm$ for a suitable prime number $p$ 
by a canonical lift of Frobenius on $E$.  So we now recall some facts
concerning canonical lifts.

\medskip

As in the introduction, let $K$ be a number field with ring of integers
$R$, and let $E/K$ be 
an elliptic curve with complex multiplication by an order $R'$ in the
imaginary quadratic field $k$.  
Let $v$ be a non-archimedean place of $K$, corresponding to a prime ideal
$\p$ of $R$ of residue characteristic $p>2$, and let $K_v$
(resp. $R_v$) denote the $v$-adic completion of $K$ (resp. $R$).  Throughout this section, 
we make the additional assumptions that $E$ has good ordinary reduction at $v$ 
and that $p$ does not divide the conductor of $R'$.  
Our hypotheses on $p$ ensure that 
there exists a canonical lift $F \in \End_{\Kbar}(E)$ of
the Frobenius endomorphism of $\Et$ (see e.g. Chapter 13, 
Theorems 12--14 of \cite{LangEF}).  The endomorphism $F$ is an isogeny
of degree $p$, and by replacing $K$ by a finite extension if 
necessary, we can and will assume that $F$ is defined over $K$.  
Also, the universal property of the N{\'e}ron model (see \cite[1.1]{Artin}) 
ensures that $F$ extends to an endomorphism $F : \E \to \E$.

Let $\E : y^2 = f(x)$ be a minimal Weierstrass model for $E/R_v$,
and let $\Et$ be the special fiber of $\E$.
We denote by 
$\Ehat$ the formal group of $\E$.

Concerning the canonical lifting of Frobenius, 
we have the following lemma:

\begin{lemma}
\label{kernelofred}
With notation and hypotheses as above (so that, in particular, $E/K$ is
an elliptic curve with complex multiplication and good ordinary
reduction at $v$),
let $I$ be the inertia subgroup of $\Gal(\Kvbar / K_v)$.  Then 
$E[F^n] \cong \mu_{p^n}$ as $I$-modules for all integers $n\geq 1$.
\end{lemma}

\pf 
We claim that $\E[F^n] = \E[p^n]^0$ as finite flat group schemes over $R$
for all $n\geq 1$.
One way to see this is as follows.  By a theorem of Deligne (see 
\cite[pp.~144-145]{Tate-FFGS}), 
$\E[F^n]$ is killed by its order $p^n$.  But since $\E[F^n]$ is connected,
it is a closed subscheme of $\E[p^n]^0$.  Since $E$ is ordinary, a consideration
of degrees forces the equality $\E[F^n] = \E[p^n]^0$, proving the claim.

By the Serre-Tate correspondence, the formal group of $\E$ is the
connected component of the $p$-divisible group of $\E$.
It follows that $\E[F^n] = \E[p^n]^0 = \Ehat[p^n]$, and the result now
follows from Lemma~\ref{formalgrouplemma}.
\qed

\medskip

We also note for future reference the fact that 
for all $P \in E(\Kbar)$, $\hhat(FP) = p\hhat(P)$.
This follows from the more general fact that for all nonzero 
$\psi \in \End(E)$, we have $\hhat(\psi(P)) = 
\deg (\psi) \hhat(P)$ for all $P \in E(\Kbar)$
(see e.g. \cite[Exercise 3.6.3]{SerreLMW}).

\medskip

As a first step toward proving Theorem~\ref{maintheorem}, 
we consider the special case where $L/K$ is an abelian extension 
which is {\em unramified} above $v$.
In this case, let $\s \in \Gal(L/K)$ be the 
Frobenius element at $v$.
(This is well-defined since $L/K$ is abelian.)

\begin{lemma}
\label{unramlemma}
Let $L/K$ be an abelian extension, unramified above $v$, and
suppose $P\in E(L)$ is not a
torsion point.  Then $FP \neq \s P$, and
\[
\l_w(FP - \s P) \geq \log p
\]
for all places $w$ of $L$ dividing $v$.  

\end{lemma}

\pf Suppose $FP = \s P$.  Then by induction, we have $F^k P = \s^k P$ for all
$k\geq 1$.  If $M$ denotes the order of $\s$ in $\Gal(L/K)$, then
we have $(F^M - 1) P = 0$.  
But $F^M - 1$ is an isogeny, which
contradicts the fact that $P$ is not a torsion point.  

For the second part of the lemma, we note that
for $Q \in E(L),  Q\neq 0$, we have $\l_w(Q) > 0 $ 
if and only if $Q$ is in the kernel of reduction mod $w$.
Also, since $L_w/K_v$ is unramified, the inequality $\l_w(Q) > 0$
implies that $\l_w(Q) \geq \log p$.
The lemma follows from this last observation, since 
$FP - \s P$ is in the kernel of
reduction mod $w$ by the definition of $F$ and $\s$.
\qed

\begin{remark}
Conversely, if $P\in E(L)$ is a torsion point 
then $FP = \s P$.  So Lemma~\ref{unramlemma} gives a
useful criterion for distinguishing non-torsion points from torsion points.
\end{remark}

We now consider the special case of Theorem~\ref{maintheorem} 
in which $K_v = \Qp$ and 
$L_w = \Q_p(\zeta_m)$ for some positive integer $m$ divisible by $p$.  
Let $\tau \in \Gal(L/K)$ be a generator of 
$\Gal(\Q_p(\zeta_m) / \Q_p(\zeta_{m/p}))$, which can be identified
with a subgroup of the inertia group $I$ of $\Gal(L_w/K_v)$.

The following elementary observation is similar to \cite[Lemma~2(1)]{AD}:

\begin{lemma}
\label{congruence}
With $\tau,m$ as above, we have
$\tau (\alpha)^p \equiv \alpha^p \bmod p$ 
for all $\alpha \in \Zp [\zeta_m]$.
\end{lemma}

\pf Write $\alpha = h(\zeta_m)$ for some polynomial $h \in \Zp [x]$.
Then 
\[
\alpha^p \equiv h(\zeta_m^p) \equiv h(\tau \zeta_m^p) \equiv
\tau \alpha^p \bmod p.
\]
\qed

The next lemma will allow us to prove Theorem~\ref{maintheorem} via
induction on the ramification of $L/K$.  It is modeled after 
Lemma~2(2) of \cite{AD}.

\begin{lemma}  
\label{ramlemma}
Let $m = p^k m'$ be an integer with $(m',p) = 1$,
let $L/K$ be a finite abelian extension, and let $w$ be a place of $L$
lying over the place $v$ of $K$.
Assume that $L_w \cong \Qp(\zeta_m)$, $K_v \cong 
\Qp$, and that $\tau$ is a generator of 
$\Gal(\Qp(\zeta_m)/\Qp(\zeta_{m/p}))$.  
Assume furthermore that $L \supseteq K(E[F^k])$.
Suppose $P \in E(L)$ is a non-torsion point, 
and let $P' = F\tau P - F P$. 
\begin{list1}
\item If $P' \neq 0$, then
$\l_w (P') \geq \log p.$
\item If $P' = 0$, then there exists a torsion point
$T \in E[F^k]$ such that $P + T$ is fixed by $\tau$.
\end{list1}
\end{lemma}

\pf For the first statement, we view $P$ as a point of $E(\Qp(\zeta_m))$.
Let $\Qpunr$ denote the completion of the maximal unramified extension
of $\Qp$, and let $I$ be the inertia group of $L_w / \Qp$.
According to Lemma~2 of \cite{Serre-Tate}, 
the reduction map gives an isomorphism
\[
E[m](\Qpunr) \cong \Et[m](\Fpbar).
\]
It follows that there exists a torsion point
$S \in E[m](\Qpunr)$ such that $P - S$ is not in the kernel of
reduction.  Since $\tau$ fixes 
torsion points unramified at $v$ and $F$ is defined over $K$, we have
\[
(\tau F - F)P = (\tau F - F)(P - S),
\]
and therefore we may assume without loss of generality for the first
part of the lemma that $P$ is not in the kernel of reduction.
Therefore we may write $P = (x,y)$ with $x,y \in \Zp[\zeta_m]$.

Since $\tau$ is in the inertia group $I$, it follows that $FP$ and 
$\tau FP$ are congruent modulo the maximal ideal $\M$ of
$\Zp[\zeta_m]$.  
We claim that they are in fact
congruent mod $\M^e = (p)$, which will establish part (1) of the
lemma.


To see this, simply note that in terms of Weierstrass coordinates, we have
$\tau F P \equiv (\tau x^p, \tau y^p)$ and $FP \equiv (x^p,y^p)$ (mod
$p$).  The claim therefore follows from Lemma~\ref{congruence}.

\medskip

For the second part of the lemma, suppose that $F\tau P = FP$.  Then
$\tau P - P \in E[F]$, which is a cyclic group of order $p$, 
so we have $\tau P - P = T'$, where $T'$ is a 
torsion point of order dividing $p$.  If
$T' = 0$ then we can take $T=0$ in the statement of the lemma, so
without loss of generality, we can suppose that $T'$ has order $p$.

By Lemma~\ref{kernelofred}, for each $n\geq 1$ we have $E[F^n] \cong 
\mu_{p^n}$ as $I$-modules.
If we fix a torsion point $U \in E[F^k] \backslash 
E[F^{k-1}]$, then it follows that 
$T'' := \tau U - U$ has order $p$, so that $T' = rT''$ for some integer
$0 < r < p$.  Define $T := -rU$.  Then
\[
\tau (P+T) - (P+T) = (\tau P - P) - r(\tau U - U) = T' - rT'' = 0,
\]
as desired.
\qed

\begin{remark}
The second part of the lemma can be interpreted in terms of Galois
cohomology.  
Let $G'$ be the cyclic group $\Gal(\Qp(\zeta_m)/\Qp(\zeta_{m/p}))$.
Since $E[F^k] \cong \mu_{p^k}$, 
the fact that for every $T' \in E[F]$
there exists $T \in E[F^k]$ such that $T' = \tau T - T$ is equivalent
to the fact that the Galois cohomology group 
$H^1(G',\mu_m)$ is killed by $m/p$.
\end{remark}

\section{Proof of Theorem~\ref{maintheorem}}
\label{proofsection}

It is clear that in the proof of Theorem~\ref{maintheorem} 
we may replace $K$ by a finite
extension, so we assume from now on that $K$ contains the field $k$ of
complex multiplication.  Also, without loss of generality we may assume
that $E$ has everywhere good reduction over $K$ 
(see \cite[Theorem~II.6.1]{Silverman2}).  This implies
(see \cite[Remark~VI.4.1.1]{Silverman2})
that for every finite extension $L/K$ and every non-archimedean place
$w$ of $L$, we have $\lambda_w(Q) \geq 0$ for all nonzero points
$Q \in E(L)$.

Let $D = [K : \Q]$.  We choose once and for all (using the Cebotarev
density theorem)
a non-archimedean place $v_0$ of $K$ of residue 
characteristic $p$ satisfying the following conditions:

\begin{list1}
\item $E$ has good reduction at $v_0$
\item $p$ splits completely in $K$
\item $p > {\rm exp}(DC_1)$, 
where $C_1$ is the constant appearing in
the statement of Corollary~\ref{Elkiescor}.  (Equivalently, this inequality
says that $A_p := (\log p)/D - C_1 > 0$.)
\end{list1}

Note that since $k\subseteq K$, hypothesis (2) automatically implies that $E$ has
ordinary reduction at $v_0$.

We will prove the following more explicit version of 
Theorem~\ref{maintheorem}:

\begin{theorem}
\label{precise}
With $K$ and $p$ as above, 
let $E/K$ be an elliptic curve having complex
multiplication.  Let $C_1$ be the constant from Corollary~\ref{Elkiescor}.
Let $L/K$ be an abelian extension, and let 
$P\in E(L)$ be a non-torsion point.  If we set 
$A_p = (\log p)/D - C_1 > 0$, then 
\[
\hhat (P) \geq \left\{
\begin{array}{ll}
\frac{A_p}{2(p+1)} &
{\rm \; if \; } L/K {\rm \; is \; unramified 
\; above \;} v_0 \\
\frac{A_p}{4p} & {\rm \; if \; } L/K {\rm \; is \; 
ramified \; above \;} v_0. \\ 
\end{array} \right.
\]
\end{theorem}


\pf 
We may assume without loss of generality that the extension $L/K$ is finite.
The proof proceeds by induction on the power of $p$ dividing the
local conductor at $v_0$ of $L/K$.  

More concretely, fix a place $w$ of $L$ lying over $v_0$, and 
let $m$ be the smallest
positive integer such that $L_w \subseteq \Q_p(\zeta_m)$.  
Such an $m$ exists
by the local Kronecker--Weber theorem, since $K_{v_0} \cong \Q_p$.
If we write $m=p^r m'$, with $(p,m')=1$, then we define the
{\em local conductor} $\f(L) = \f_{v_0}(L/K)$ to be $p^r$.  This coincides
with the definition of the conductor from local class field theory.
As the notation suggests, $\f_{v_0}(L/K)$ is independent of the 
chosen place $w \mid v_0$.

\medskip

Let $P \in E(L)$ be a non-torsion point.
If $\f(L)=1$ (i.e., $L/K$ is unramified above $v_0$), then 
by Lemma~\ref{unramlemma}
(letting $\s$ be the Frobenius element of $L/K$ above $v_0$) we have
\[
\l_w (\s P - FP) \geq \log p
\]
for all places $w$ of $L$ with $w | v_0$.

It follows that
\[
\sum_{w | p} d_w \l_w (\s P - FP) \geq 
\sum_{w | v_0} d_w \log p = d_{v_0} \log p = \frac{\log p}{D}.
\]

As noted above, for every nonarchimedean place $w \in M_L^0$ and
every nonzero point $Q \in E(L)$ we have
\[
d_w \l_w(Q) \geq 0.
\]

By Corollary~\ref{Elkiescor},
we also know that for every archimedean place $w \in M_L^\infty$ we have
\[
\l_w(Q) \geq -C_1.
\]

Applying these inequalities to $Q = \s P - FP$, and using 
the decomposition
\[
\hhat(Q) = \sum_{v \in M_\Q} \sum_{w | v} d_w \l_w(Q)
\]
of the global canonical height as a sum of N{\'e}ron local heights,
we find that
\[
\hhat(\s P - FP) \geq \frac{\log p}{D} - C_1 = A_p > 0.
\]

Applying the parallelogram law
\[
\hhat(P_1 - P_2) = 2\hhat(P_1) + 2\hhat(P_2) - \hhat(P_1 + P_2)
\]
and noting that $\hhat(P_1 + P_2) \geq 0$, $\hhat(\s P) = \hhat(P)$,
and $\hhat(FP) = p\hhat(P)$, we conclude that
\[
\hhat(P) \geq 
\frac{A_p}{2(p+1)}.
\]

This establishes the base case $\f(L) = 1$ of the induction.

\medskip

Assume now that $\f(L) \geq p$ (so that $L/K$ is ramified above $v_0$),
and that the theorem is true for all abelian extensions $L'$ of $K$ 
such that $\f(L') < \f(L)$.

As the local conductor of 
$L(\zeta_{m})$ over $K$ is still $\f = p^k$, we may assume without loss
of generality that $L$ contains the $m$th roots of unity, so that
for all places $w$ of $L$ lying over $v_0$ we have $L_w \cong \Q_p(\zeta_m)$.

By Lemma~\ref{kernelofred}, all points in $E[F^k]$
are defined over $\Qp(\zeta_{p^k}) \subseteq \Qp(\zeta_m) = L_w$.  
They are also defined over an abelian extension of $K$.
Without loss of generality, we may therefore also assume that 
$K(E[F^k]) \subseteq L$.

Let $\tau \in \Gal(L/K)$ be as in Lemma~\ref{ramlemma}, 
and let $P' = \tau FP - FP$.
If $P' \neq 0$, then by part (1) of Lemma~\ref{ramlemma} we have
\[
\l_w (P') \geq \log p
\]
for all places $w$ of $L$ with $w | v_0$.

An argument parallel to that in the case $\f(L) = 1$ now shows that
\[
\hhat(P) \geq \frac{A_p}{4p}.
\]

Finally, if $P' = 0$, then part (2) of Lemma~\ref{ramlemma} tells us that
there exists a torsion point $T \in E[F^k] \subseteq E(L)$ such that 
$P + T$ is defined over the fixed field $L'$ of $\tau$.  Note that
$L'_w = \Qp(\zeta_{m/p})$, and therefore $\f(L') < \f(L)$.  
Since $\hhat(P) = \hhat(P + T)$, we can apply the inductive hypothesis
to $P+T$ to conclude.
\qed


\section{Proof of Theorem~\ref{maintheorem2}}
\label{proofsection2}

Suppose $E/K$ is an elliptic curve whose $j$-invariant is non-integral
at the finite place $v_0$ of $K$.
Without loss of generality, we may replace $K$ by a finite extension and 
we henceforth assume that $E/K$ is semistable, with split multiplicative
reduction at $v_0$.
Let $\nu_0 = \ord_{v_0}(\jinv)$.

\medskip

Suppose, for the sake of contradiction, that $\{ P_i \}$ is an 
infinite sequence of points in $E(\Kab)$ such that $\hhat(P_i) \to 0$.  
Let $L_i$ be the Galois closure of $K(P_i) / K$, which by assumption is
an abelian extension of $K$.
Let $e_i$ be the ramification index above $v_0$ of $L_i/K$.

We consider two cases.  

\medskip
\underline{{\bf Case 1:}} The sequence $e_i$ is bounded above. 
\smallskip

For notational convenience in the sequel, let $e$ be an upper bound for
the integers $e_i$, and define 
$c := \frac{1}{e^2\nu_0^2} > 0$.
Pick any point $Q$ from the sequence $\{ P_i \}$, let $L$ be the
Galois closure of $K(Q)/K$,
and let $Q_1,\ldots,Q_N$ be the
Galois conjugates of $Q$ over $K$.  Note that Northcott's theorem implies that
we can choose $Q$ so that $N$ is as large as we please.

By repeatedly applying the parallelogram law
\[
\hhat(Q_i + Q_j) + \hhat(Q_i - Q_j) = 2\hhat(Q_i) + 2\hhat(Q_j)
\]
and noting that $\hhat(Q_i) = \hhat(Q_j) = \hhat(Q)$ and 
$\hhat(Q_i + Q_j) \geq 0$ for all $i,j$, 
we obtain (following \cite{HS2}) the inequality
\begin{equation}
\label{proofeqn0}
\hhat(Q) \geq \frac{1}{4N(N-1)} \sum_{i\neq j} \hhat(Q_i - Q_j).
\end{equation}

In order to make use of this observation, we need to estimate the
contribution at places lying over $v_0$ to the local decomposition
of the global canonical height as a sum of N{\'e}ron local heights.

\medskip

To do this, let $w$ be any place of $L$ lying over $v_0$.
Then $\log |j_E|_w > 0$ and $\nu := \ord_w(\jinv) = e(w/v_0)\nu_0 \leq e\nu_0$.
By Lemma~\ref{HSlemma}, 
we have, for each integer $N\geq 2$, the inequalities
\begin{equation}
\begin{array}{lll}
\label{proofeqn1}
\sum_{i\neq j} \lambda_w(Q_i - Q_j) & \geq &
\frac{1}{12} \log |j_E|_w 
\left( \left(\frac{N}{\nu}\right)^2 - N \right) \\
& \geq & \frac{c}{12}N^2\log |j_E|_w - \frac{N}{12}\log |j_E|_w, \\
\end{array}
\end{equation}
and therefore
\begin{equation}
\label{proofeqn1bis}
\frac{1}{N(N-1)} \sum_{i\neq j} \lambda_w(Q_i - Q_j) \geq
\frac{c}{12} \log |j_E|_w - \frac{1}{12(N-1)}\log |j_E|_w.
\end{equation}

\medskip

We have been working at the specific place $v_0$, but 
by the same reasoning, we see that for any place 
$w$ of $L$ where $E$ has bad (and therefore split multiplicative) reduction
we have
\begin{equation}
\label{proofeqn2}
\frac{1}{N(N-1)} \sum_{i \neq j} \l_w(Q_i - Q_j)
\geq - \frac{1}{12(N-1)}\log |j_E|_w.
\end{equation}

For any place $w$ where $E$ has good reduction, we have
\begin{equation}
\label{proofeqn3}
\frac{1}{N(N-1)} \sum_{i \neq j} \l_w(Q_i - Q_j) \geq 0.
\end{equation}

Finally, for any infinite place $w$ of $L$ we have, by 
Lemma~\ref{Elkies},
\begin{equation}
\label{proofeqn4}
\frac{1}{N(N-1)} \sum_{i \neq j} \l_w(Q_i - Q_j)
\geq - \frac{1}{12(N-1)}\log^+ |j_E|_w + o(1),
\end{equation}
where $o(1)$ denotes a term which tends to 0 as $N$ tends to infinity.

Combining inequalities (\ref{proofeqn1}) through (\ref{proofeqn4}), 
and using the decomposition of
the global height as a sum of local heights, we obtain, using Lemma~\ref{absolutevaluelemma}:

\begin{equation}
\begin{array}{lll}
\label{proofeqn5}

\frac{1}{N(N-1)} \sum_{i\neq j} \hhat(Q_i - Q_j) & = &
\sum_{w \in M_L} d_w \frac{1}{N(N-1)} 
\left(\sum_{i\neq j} \l_w(Q_i - Q_j) \right) \\
& & \geq \frac{c}{12} \sum_{w\mid v_0} d_w \log |j_E|_w
- \frac{1}{12(N-1)} \sum_{w \in M_L} d_w \log^+ |j_E|_w \\
& = & \frac{c}{12} d_{v_0} \log |j_E|_{v_0}  - \frac{1}{12(N-1)} h(j_E) \\
& = & \frac{c}{12} d_{v_0} \log |j_E|_{v_0}  + o(1). \\
\end{array}
\end{equation}

Combining this with inequality (\ref{proofeqn0}), we find that there
are positive real constants $c' = d_{v_0} \log |j_E|_{v_0} / 48$ and
$c'' = cc'$ depending only on $E,K,v_0$, and $e$ such that

\begin{equation}
\label{proofeqn6}
\hhat(Q) \geq c'' + o(1).
\end{equation}

Since $d_i \to\infty$, $Q$ can be chosen so that $N$ is as large as
we like.  However, for $N$ sufficiently large, the inequality 
(\ref{proofeqn6}) contradicts the fact that $h(P_i) \to 0$.  We therefore
conclude that Case 1 cannot occur.

\medskip

\underline{{\bf Case 2:}} $e_i \to \infty$.

\smallskip

As before, we pick an element $Q \in E(L)$ from the sequence 
$\{ P_i \}$.  Let $L$ be the Galois closure of $K(Q)/K$; since $L/K$ is
abelian, the inertia groups of $L/K$ at each place lying over $v_0$ are
all the same group $I$.  

Let $Q_1,\ldots,Q_N$ be the conjugates of $Q$ under $I$.  Note that
the hypothesis $e_i \to \infty$ guarantees that we can choose $Q$ so
that $N$ is as large as we like.

Since the $Q_i$'s are conjugate under the inertia group $I$, we see that
for every place $w$ of $L$ dividing $v_0$, 
$Q_i - Q_j$ is in the kernel of reduction mod $w$ for all $i,j$.  It
follows that 
for all $i\neq j$ we have
\[
\lambda_w(Q_i - Q_j) \geq 
\frac{1}{12}\log |j_E|_w
\]
and therefore
\[
\frac{1}{N(N-1)} \sum_{i\neq j} \lambda_w(Q_i - Q_j) \geq
\frac{1}{12}\log |j_E|_w
\]
for all $w \mid v_0$.

Using the parallelogram law and Lemma~\ref{absolutevaluelemma} as before, we find that
\[
\begin{array}{lll}
\hhat(Q) & \geq & \frac{1}{4}\sum_{w|v_0}
d_w (\frac{1}{12}\log |j_E|_w) + o(1) \\
& = & 
\frac{1}{48} d_{v_0} \log |j_E|_{v_0} + o(1) \\
& = & c' + o(1), \\
\end{array}
\]
which for $N$ large contradicts the fact that $\hhat(P_i) \to 0$.
\qed


\section*{Acknowledgements}

The author would like to thank Brian Conrad for helpful discussions
and Robert Rumely for his useful comments on an earlier version of
this manuscript.


\end{document}